\documentclass[12pt]{amsart}
\usepackage{amscd,amsthm,amsmath} % SRC Specials: DVI [Inverse] Search
\usepackage{pb-diagram}
\newtheorem{thm}{Theorem}[section]
\newtheorem{cor}[thm]{Corollary}
\newtheorem{lem}[thm]{Lemma}
\newtheorem{prop}[thm]{Proposition}
\newtheorem{examp}[thm]{Example}
\newtheorem{ques}[thm]{Question}
\newtheorem{conj}[thm]{Conjecture}
\theoremstyle{definition}
\newtheorem{defn}[thm]{Definition}
\theoremstyle{remark}
\newtheorem{rem}[thm]{Remark}
\numberwithin{equation}{section}
\begin{document}

\title[Forgetable map]{Phantom maps and Forgetable maps}%
\author{Jianzhong Pan }%
\address{Institute of Math.,Academia Sinica ,Beijing China and
 Department of Mathematics Education , Korea University , Seoul , Korea }%
\email{pjz@kuccnx.korea.ac.kr}%
\author{Moo Ha Woo}
\address{ Department of Mathematics Education , Korea University , Seoul , Korea}
\thanks{The first author is partially supported by the NSFC project
19701032 and ZD9603 of Chinese Academy of Science and Brain
Pool program of KOSEF and the second author wishes to acknowledge
the financial supports of the Korea Research Foundation made in the
program year of (1998) and TGRC 99}%

\date{July 13,1999}%
%\dedicatory{}%
%\commby{}%
% ----------------------------------------------------------------
\begin{abstract}
 In this note, we attack   a question posed
ten years ago by Tsukiyama about the injectivity of the so- called
Forgetable map. We show that we can insert the Forgetable map in
an exact sequence and that  the problem can be reduced to the
computation of the sequence which turns out  unexpectedly to be
related to the phantom map problem and the famous Halperin
conjecture in rational homotopy theory.
\end{abstract}
\maketitle
% ----------------------------------------------------------------
 \section{Introduction} \label{S:intro}
Homotopy theory is now an area which contains many research
topics. Among them are two interesting ones : Phantom map theory
which is well known and Forgetable map which is maybe not so well
known. An interesting relation have been observed in this paper
which is used to attack the Forgetable map problem posed by
Tsukiyama\cite{kt1},\cite{kt2} ten years ago.

A pair of  maps $f$ and $g$ from  a CW complex $X$ to a
topological space $Y$ is called a phantom pair if the restriction
of $f$ and $g$ to the n-skeleton of the complex $X$ is homotopic
for all $n \geq 0$.In this case , we call the map $f$ a $g$-
phantom map with respect to the map $g$. Denote by $Ph^g(X,Y)$ the
set of homotopy classes of $g$-phantom maps from $X$ to $Y$. It is
clear that the concept of $g$-phantom map is  homotopy invariant.
Especially , if $g$ = constant map, then the $g$-phantom mapwas
called just a phantom map and if  $g=id$, then $f$ is a
$g$-phantom map if and only if $f$ is a weak identity as defined
by Roitberg\cite{roit1}.

Historically Adams and Walker\cite{aw} found the first nontrivial
phantom map and Gray made the first detailed study of phantom maps
in his PhD thesis\cite{gray} . Later many other authors ,
McGibbon, Meier , Moller , Oda-Shitanda , Roitberg, Sullivan,
Zabrodsky , etc. contributed a lot of ideas to this beautiful
area,see \cite{cm95} for a comprehensive survey of this area.
Among many others,they proved the following   which are crucial to
our later application .
\begin{thm}\cite{cm95}
Let X and Y be nilpotent CW complexes with finite type. Then the
set $Ph^g(X,Y)$ is either one point set or uncountable.
\end{thm}
In general,$Ph^*(X,Y)$ is a proper subset of $[X,Y]$ . But the
following theorem says they are equal under some conditions and
even more is true.
\begin{thm}\label{T:thm1.4} \cite{cm95}
Let $X$ and $Y$ be nilpotent CW complexes with finite type . If
$[X_{\tau},Y]=[\Sigma X_{\tau},Y]=*$, then
\begin{equation}
Ph^*(X,Y)=[X,Y]=[X_{(0)},Y]=\prod_{n >0}
H^n(X,\pi_{n+1}(Y)\bigotimes R) \label{eq:phan}
\end{equation}
\end{thm}
\begin{thm} \label{T:thm1.5} \cite{cm95}
Let $Y=\Omega^nK$ and $X=\Sigma^mZ$ for $m , n\geq 0$ where $K$ is
a simply connected finite CW complex  .If $Z$ satisfies one of the
following conditions, then we have
\begin{equation}
Ph^*(X,Y)=[X,Y]=[X_{(0)},Y]=\prod_{n>0}H^n(X,\pi_{n+1}(Y)\bigotimes
R)
\end{equation}
\begin{itemize}
\item{$Z=BG$,$G$ is a connect compact lie group with finite
fundamental group ,or}
\item{$Z$ is a connected infinite loop
space with finite fundamental group ,or}
\item{$Z$ is a
1-connected finite Postnikov space, i.e. $\pi_iZ=0$ if $i$ is
large enough}.
\end{itemize}
\end{thm}
In the above theorem the target space is the (iterated) loop space
of a finite CW complex. To deal with some essential infinite space
Zabrodsky extend the above theorem as follows
\begin{thm}\label{T:zab} \cite{az}
The equation (\ref{eq:phan}) remains true if $X=\Sigma^mK(H,l+2)$
and $Y=\Omega^nBaut(K)$  where $m,l,n\geq 0$, $H$ is an abelian
group and K is a finite CW complex.
\end{thm}

Theorem \ref{T:thm1.5} and \ref{T:zab} say, in some case ,all maps
are phantom maps and the homotopy classes of them can be
calculated . The general phantom pair is studied only briefly by
Oda-Shitanda\cite{os} and seems to be forgot later.
Roitberg\cite{roit1,roit2} has studied the weak identities and
posed several interesting open questions about them and later
Shitanda had also some related works on it.

According to our point of view the main problem one faces with the
general phantom pair is the following
\begin{ques}\label{T:ques1}
Let  $g_1,g_2:X \to Y $ are two maps. What is the implication
relation between $Ph^{g_1}(X,Y)=\{g_1\}$ and
$Ph^{g_2}(X,Y)=\{g_2\}$
\end{ques}
A well known result in this direction is
\begin{thm} \label{T:homo}\cite{os}
If Y is a H-space with inverse,then for any two maps $g_1,g_2:X
\to Y $ the equations in Question\ref{T:ques1} are equivalent.
\end{thm}
In our application we have to extend the Theorem \ref{T:homo} to
the case where $Y$ is not a H-space. Actually what we need is
about somewhat more general notion.See
Theorem\ref{T:section2maincoro1},\ref{T:liephan} for details.

 Now let us turn to  Forgetable map. Given
a principal G-bundle $\pi:P \to B$, Let
\[
aut^G(P)=\{g | g:P \to P  \text{  is a G-equivariant homotopy
equivalence}\}
\]
and
\[
aut(P)=\{g | g:P \to P  \text{  is a homotopy equivalence}\}
\]
There is a natural map $f:aut^G(P)\to aut(P)$. Let
\[
Aut^G(P)= \pi_0(aut^G(P))
\]
  and
\[
  Aut(P)= \pi_0(aut(P))
\]
Then the map $f$ induces   a map
\[
F:Aut^G(P) \to Aut(P)
\]
which is called a Forgetable map by Tsukiyama. The question posed
by Tsukiyama in \cite{dwk} is the following
\begin{ques} \label{T:ques2}
Is the forgetting  map $F$   injective?
\end{ques}
In  \cite{kt1,kt2},Tsukiyama constructed examples which answers
negatively the question \ref{T:ques2}  and gave a sufficient
condition which answers positively the question \ref{T:ques2} .
His example is the following:

\begin{examp}
Given a connected compact Lie group, which is not a torus, $G$ and
the maximal torus $T$. There is principal $G$-bundle $G \to G/T
\to BT$ over $BT$  which is classified by the natural map $Bi:BT
\to BG$ where $i:T \to G$ is the inclusion of $T$ into $G$. Then
$Aut(G/T)$ is finite and there is an exact sequence
\[
0 \to \pi_1(map(BT,BG),Bi) \to Aut^G(G/T) \to Aut(BT)
\]
Since $\pi_1(map(BT,BG),Bi)$ is uncountable, $ Aut^G(G/T)$ is
uncountable and   thus the Forgetable map $F:Aut^G(G/T) \to
Aut(G/T)$ is not injective.
\end{examp}
 One of the main results in this paper
is the following
\begin{thm}\label{T:main5}
Let $\pi:P \to:B$ be a principal $G$-bundle with $P$ a finite CW
complex.Then there is an exact sequence
\[
\pi_1aut(P) \to \pi_1(map_*(BG,Baut(P)),c) \to Aut^G(P)
\overset{F}{\to} Aut(P)
\]
where  $c:BG \to Baut(P)$ is determined by the principal bundle.
\end{thm}
\begin{rem}
In the above theorem , the calculation of the $KerF$ is in some
sense equivalent to the calculation of  the group
\[
\overline{G}=: \pi_1(map_*(BG,Baut(P)),c)
\]

\begin{itemize}
\item{$\overline{G}=0 \Rightarrow KerF=0$}
\item{$\overline{G}$ is uncountable $\Rightarrow KerF$ is uncountable
since
\[
\pi_1aut(P)=\pi_1(map(P,P),id)
\]
 is countable.}
\end{itemize}
\end{rem}

The following Theorem shows that, in certain case , $KerF$ is
either zero or uncountable .
\begin{thm}\label{T:4}
Let $\pi:P \to B$ be  as above , $k:B \to BG$ the classifying map,
$\bar{k}:\bar{B} \to BG$ the associated fibration with fiber
$\bar{P}$ and $c:BG \to Baut(P)$ the classifying map of $\bar{k}$.
The followings are true.
\begin{itemize}
 \item{If $KerF$ is uncountable if $Ph_1^c(BG,Baut(P))\neq0$}
 \item{If $c$ is phantom map and $G$ be a connected compact Lie group , $KerF$ is either zero or uncountable}
\item{If $G$ be a 1-connected K(H,m),$KerF$ is either zero or uncountable}
\end{itemize}
\end{thm}

According to Tsukiyama\cite{kt1}\cite{kt2}, it is possible that
$KerF$ is zero or uncountable.The above theorem says in certain
case this is the only possibilities. Furthermore we will show that
results in phantom map theory  and  rational calculation which is
usually not so difficult  can be used to decide when the $KerF$ is
zero or uncountable. The  theorem above leads to  another natural
question
\begin{ques}
Is it possible that the $KerF$ is finite or countable?
\end{ques}
 Now we will give concrete conditions for injectivity or noninjectivity of the
 Forgetable map. First we assume the group $G$ to be a connected
 compact Lie group.
\begin{thm}\label{T:5}
Let $P$ be a 1-connected CW complex  . Then the followings are
equivalent:
\begin{itemize}
\item{there is a a connected compact Lie group $G$ and a principal
$G$-bundle such that the total space has the homotopy type of $P$
, the classifying map c is a phantom map and the associated
Forgetable map $F$ has uncountable kernel}
 \item{$\bigoplus_{i>0}\pi_{2i+1}(map(P_{(0)},P_{(0)});id)$ is nontrivial}
\end{itemize}
\end{thm}
Before giving concrete examples of principal bundles with
noninjective Forgetable map,we recall some backgrounds. A space
$P$ is called elliptic if $H^i(X,Q)=0$   and $\pi_i(X)\bigotimes
Q=0$ when $i$ is sufficient large. Let $(X,*)$ be any pointed
space. The Gottlieb group (or evaluation subgroup)\cite{gottlieb1}
is defined by
\[
G_n(X)=Im\{ev_*:\pi_n(map(X,X);id) \to \pi_n(X)\}
\]
where $ev:map(X,X)_{id}  \to  X$ is defined by $ev(f)=f(*)$

The Gottlieb groups are extremely difficult to compute in general.
However for rational space there have been some remarkable results
on the Gottlieb groups.
\begin{thm}\cite{felix}
If $X$ is finite CW complex,then $G_{even}(X_{(0)})=0$
\end{thm}
\begin{thm}\cite{smith}
Let $X$ be a finite 1-connected complex,then $G_*(X_{(0)})  \neq
0$ if $X$ is rationally nontrivial and elliptic.
\end{thm}
Now the following Theorem follows immediately from the last three
theorems.
\begin{thm}
Let $P$ be any  1-connected elliptic finite CW complex.Then there
is a compact Lie group $G$ and  a principal $G$-bundle $\pi:P' \to
B$ such that the Forgetable map has uncountable kernel and $P,P'$
has the same homotopy type.
\end{thm}
\begin{rem}
We don't know if there exists  1-connected CW complex $P$ such
that $\bigoplus_{i>0}\pi_{2i+1}(map(P_{(0)},P_{(0)});id)=0$.
\end{rem}
If we assume $G=K(H,2m)$ where $H$ is a finitely generated abelian
group and $m \geq 1$.We have the following
\begin{thm}\label{T:1}
Let $P$ be a 1-connected CW complex . Then the followings are
equivalent:
\begin{itemize}
\item{for all $m \geq 1$ , finitely generated abelian group $H$ and
every principal $K(H,2m)$-bundle with total space homotopy
equivalent to $P$ ,the associated Forgetable map is injective}
\item{$ \bigoplus_{i>1}\pi_{2i}(map(P_{(0)},P_{(0)});id)=0$}
\end{itemize}
\end{thm}
Now we want to give some examples with injective Forgetable map.
Again we first recall some backgrounds. For simplicity we will
assume spaces involved are 1-connected.

A space $X$ is said to be of type $F_0$ if $\dim H^*(X;Q)<
\infty$, $\dim \pi_*(X)\bigotimes Q< \infty$ and $H^{odd}(X;Q)=0$.
One of the most beautiful conjectures in rational homotopy theory
is the following
\begin{conj}\cite{dupont}
Let $P \to E \to B$ be a fibration such that the fiber $P$ is
homotopy equivalent to a CW complex of type $F_0$. Then the Serre
spectral sequence(with rational coefficients in Q) of the
fibration collapses at the $E^2$ term.
\end{conj}
In his two remarkable papers \cite{meier1,meier2}, W.Meier found
the relation between Halperin conjecture and the vanishing of the
$\pi_{even}map((P_{(0)},P_{(0)});id)$ , i.e.
\begin{thm}\label{T:2}
Let $P$ be of type $F_0$. Then the followings are equivalent.
\begin{itemize}
\item{The Serre
spectral sequence of every fibration with fiber $P$ collapses at
the $E^2$ term}
\item{$\pi_{even}map((P_{(0)},P_{(0)});id)=0$}
\end{itemize}
\end{thm}
The Halperin conjecture have been verified for a number of special
cases. The result obtained so far can be stated as follows

\begin{thm}\label{T:3}\cite{bl,lupt,shte,thomas}
Let $P$ be a space satisfying one of the following conditions.
Then the Halperin conjecture is true:
\begin{itemize}
\item{$P$ is a Kahler manifold}
\item{$H^*(P;Q)$ as an algebra has at most 3 generators}
\item{$P=G/U$ where $G$ is a compact Lie group and $U$ is a closed subgroup of maximal rank}
\end{itemize}
\end{thm}
Comparing Theorem\ref{T:1} and Theorem\ref{T:2} we obtain
immediately the following
\begin{thm}
Let $P$ be a 1-connected finite CW complex of type $F_0$. Then the
followings are equivalent:
\begin{itemize}
\item{The Halperin conjecture is true for $P$}
\item{For every $m \geq 1$ ,finitely generated abelian group $H$ and every principal
$K(H,2m)$-bundle with total space homotopy equivalent to $P$ ,the
associated Forgetable map is injective and
$\pi_2map((P_{(0)},P_{(0)});id)=0$}
\end{itemize}
\end{thm}
\begin{cor}
Let $P$ be a 1-connected finite CW complex satisfying the
condition of Theorem\ref{T:3}. Then For all $m \geq 1$ ,finitely
generated abelian group $H$ and every principal $K(H,2m)$-bundle
with total space homotopy equivalent to $P$ ,the associated
Forgetable map is injective
\end{cor}

 In section 2 we will introduce the
phantom element which is a generalization of phantom pair and use
this concept to prove  Theorem\ref{T:section2maincoro1} and
Theorem\ref{T:liephan}. In section 3, we will study the forgetting
map and try to insert it into an exact sequence and prove the
Theorem\ref{T:main5}. In section 4, the results of previous
sections will be applied here to get the precise information about
the Forgetable map and prove Theorem\ref{T:4}, \ref{T:5},
\ref{T:1}. In this paper all our basic spaces will be assumed to
be CW complexes with finite type . We will also use the following
notations:
\begin{itemize}
\item{$X_n$ is the $n$-th skeleton of X}
\item{$map(X,Y)$ is the space of continuous mappings from $X$ to $Y$}
\item{$map_*(X,Y)$ is the  subspace of pointed mappings from $(X,x_0)$ to $(Y,y_0)$}
\item{$l:X \to X_{(0)}$ is the rationalization}
\item{ Let $\tau:X_{\tau} \to X$ be the homotopy fiber of $l$.
Then $X_{\tau} \overset{\tau}{\to} X \to X_{(0)}$ is a cofibration
up to homotopy }
\item{$\hat{e} : Y \to \hat{Y}$ is the profinite completion of Sullivan \cite{sullivan}}
\end{itemize}
The first author thanks the Brain Pool Program of KOSEF for the
support to his visit to Korea University and the Prof.Woo and  the
Department of Mathematics Education ,Korea University for the
hospitality during which this work was completed.Thanks is also
due to Prof.Shen Xinyao who brought me to the area of topology ten
years ago and Prof.McGibbon who kindly sent me his papers and
preprints.

\section{Phantom element}
Let us begin with some definitions Let $X$ be a CW complex,$Y$ a
space and  $f,g:X \to Y$ two maps .A map $f:X \to Y$ is called a
phantom map if $f|_{X^n}$(the restriction of $f$ to the n-skeleton
of $X$) is homotopic to the constant map for all $n \geq 0$. In
\cite{os} , Oda and Shitanda defined that  $f $ and $g $ are a
phantom pair if  $f|_{X^n}$ and  $g|_{X^n}$ are homotopic for all
$n \geq 0 $. For a fixed map $g:X \to Y $ we denote by $Ph^g(X,Y)$
the set of homotopy classes of maps $f$ such that  $f $ and  $g $
are phantom pair. Each element of $Ph^g(X,Y)$  is also called a
$g$-phantom map.

Here we generalize the concept of the phantom pair as follows
\begin{defn}
Let spaces $X$ be a CW complex,$Y$ be a space  and $g:X \to Y$ any
map. Then an element $\alpha \in\pi_j(map_*(X,Y);g) $ is called a
$g$-phantom element if $(i_n^*)_{\#}(\alpha)=0$ for all $n \geq 0$
where $(i_n^*)_{\#}: map_*(X,Y)\to map_*(X^n,Y)$ is the
homomorphism induced by the inclusion $i_n :X^n \to X$. Denoted by
\[
Ph_j^g(X,Y)=\{\alpha \in\pi_j(map_*(X,Y);g)| \alpha \text{ is a
$g$-phantom element}  \}
\]
\end{defn}

If  $j=0$,then $\alpha$ is a $g$-phantom element iff it represents
the homotopy class of a map which is a $g$-phantom map. If
$g$=constant map,then a $g$-phantom map is the same as phantom
map. Since Adams and Walker\cite{aw} found the first essential
phantom map, this area has attracted interests of many
mathematicians.

Let us recall some basic results about homotopy of a sequence of
fibrations at first. Let
\[
\cdots \to X_n \overset{\pi_n}{\to} X_{n-1} \to \cdots
\]
be a sequence of fibrations of spaces and $X $ be the inverse
limit of the above inverse system. If we choose base points $x_n
\in X_n$ such that $\pi_n(x_n)=x_{n-1}$.It was shown by Bousfield
and Kan in \cite{bk} that there exists the following short exact
sequence for $j \geq 0$
\[
* \to \underset{\leftarrow}{\lim}^1_n\pi_{j+1}(X_n,x_n) \to
\pi_{j}(\underset{\leftarrow}{\lim}X_n,x_n) \to
\underset{\leftarrow}{\lim}\pi_j(X_n,x_n) \to *
\]
Now if $X$ is a CW complex with skeleton $X^n$ and $Y$  a space
,then
\[
\cdots \to map_*(X^n,Y) \overset{i_{n,n-1}^*}{\to}map_*(X^{n-1},Y)
\to \cdots
\]
is a sequence of fibrations with $map_*(X,Y)$ as the inverse
limit.

\begin{cor}
Let $X$ and $Y$ be nilpotent CW complexes of finite type and $g:X
\to Y $ be any map. Then for all $j \geq 0$ there exists a short
exact sequence
\[
* \to \underset{\leftarrow}{\lim}^1_n\pi_{j+1}(map_*(X^n,Y);g|X^n)
\to \pi_j(map_*(X,Y);g) \to
\]
\[
\underset{\leftarrow}{\lim}_n\pi_j(map_*(X^n,Y);g|X^n) \to *
\]
\end{cor}

\begin{cor} \label{T:descofph}
Let $X$ and $Y$ be nilpotent CW complexes of finite type and $g:X
\to Y $ be any map. Then for all $j \geq 0$,we have
\[
Ph^g_j(X,Y)=
\underset{\leftarrow}{\lim}^1_n\pi_{j+1}(map_*(X^n,Y);g|X^n)
\]
If $j=0$, we recover
$Ph^g(X,Y)=\underset{\leftarrow}{\lim}^1_n\pi_1(map_*(X^n,Y);g|X^n)$
\end{cor}

A natural problem about $Ph^g_j(X,Y)$ is  its cardinality . For
this we have the following

\begin{thm}\cite{cm95}
The first derived inverse limit of an inverse system of countable
groups  is either one point set or uncountable.
\end{thm}\label{T:thm2.7}
\begin{cor}
Let  $X$ and $Y$ be nilpotent CW complexes of finite type and $g:X
\to Y $ be any map. Then $Ph_j^g(X,Y)$ is either one point set or
uncountable for all $j \geq 0$.
\end{cor}

Another natural question is the extended version of Question
\ref{T:ques1}.
\begin{ques}
For two maps $f,g:X \to Y$,what is the relation between
$Ph_j^g(X,Y)$ and $Ph_j^f(X,Y)$?
\end{ques}
 The first result in this direction is an extension of the result
of Oda-Shitanda\cite{os}.
\begin{thm}
Let X be a CW complex and Y be an H-space with inverse. Then we
have $Ph_j^g(X,Y) = Ph_j^f(X,Y)$
\end{thm}
\begin{proof}
If $j=0$, this is Theorem 3.7(2) of Oda-Shitanda\cite{os} and if
$j>0$, this follows from the following lemma and the Corollary
\ref{T:descofph}.
\end{proof}
\begin{lem}
Let $i:X_1 \to X_2$ be a map and Y be an H-space with inverse.
Then for any map $f:X_2 \to Y$ the following diagram is
commutative up to homotopy
\[
\begin{CD}
map_*(X_2,Y)_f @>>i^*> map_*(X_1,Y)_{f\circ i}\\
 @VVh_2V      @ VVh_1V \\
 map_*(X_2,Y)_* @>>i^*> map_*(X_1,Y)_*
\end{CD}
\]
where $h_1, h_2$ are defined during the proof.
\end{lem}

\begin{proof}
Since $Y$ is an H-space  with inverse, there is a  multiplication
$\mu:Y \times Y \to Y$ , a map  $ \nu:Y \to Y$ is called an
inverse for $\mu$ if each composite map
\[
Y \overset{(1, \nu)}{\to } Y \times Y \overset{\mu}{\to}Y \text{
and }Y \overset{( \nu,1)}{\to } Y \times Y \overset{\mu}{\to}Y
\]
are homotopic to the constant map $*:Y \to Y$. From these we can
construct the maps $h_1,h_2$ by the composite
\[
map_*(X_j,Y)_{f_j}\overset{(1, \nu_*)}{\to} map_*(X_j,Y)_{(f_j}
\times map_*(X_j,Y)_{ \nu \circ f_j}\to
\]
\[
 \to map_*(X_j,Y\times Y)_{(f_j, \nu \circ f_j)}
\overset{\mu_*}{\to}map_*(X_j,Y)_*
\]
where $f_2=f$ and $f_1=f\circ i$

It is easy to see from the definition of the maps $h_1,h_2$ that
the diagram in the lemma   is commutative up to homotopy.
\end{proof}

\begin{thm}\label{T:fundamental}\cite{pw}
Let $X,Y$ be nilpotent CW complexes of finite type and $f:X \to Y
$ be any map.Then $\alpha \in \pi_j(map_*(X,Y);f)$ is a phantom
element $\Leftrightarrow$
\begin{itemize}
\item{$(\hat{e}_*)_{\#}(\alpha)=0$ where $(\hat{e}_*)_{\#} : \pi_j(map_*(X,Y);f) \to  \pi_j(map_*(X,\hat{Y});\hat{f})$}
\item{$\tau^*_{\#}(\alpha)=0$  where $\tau^*_{\#}:\pi_j(map_*(X,Y);f) \to  \pi_j(map_*(X_{\tau},Y);f_{\tau})$}
\end{itemize}
\end{thm}

\begin{proof}
We will only prove the "if" part which is necessary for our
application in this paper. For the proof of the other parts , see
\cite{pw}.

Let  $\alpha \in \pi_j(map_*(X,Y);f)$ such that
$(\hat{e}_*)_{\#}(\alpha)=0$. If we consider the following
commutative diagram
\[
\begin{CD}
\pi_j(map_*(X,Y);f) @>>(i_n^*)_{\#}> \pi_j(map_*(X^n,Y);f|_{X^n})
\\
 @VV(\hat{e}_*)_{\#}V  @VV(\hat{e}_*)_{\#}V  \\
\pi_j(map_*(X,\hat{Y});\hat{f}) @>>(i_n^*)_{\#}>
\pi_j(map_*(X^n,\hat{Y});\hat{f}|_{X^n})
\\
\end{CD}
\]
then we have  $(\hat{e}_*)_{\#}\circ (i_n^*)_{\#}(\alpha)=
(i_n^*)_{\#}\circ (\hat{e}_*)_{\#}(\alpha)=0$. In \cite{sullivan},
Sullivan showed that if $Y$ is a nilpotent space , $\hat{e}:Y \to
\hat{Y}$ and $h,g:Z \to Y$ be any two maps where $Z$ any finite CW
complex such that $\hat{e}\circ g\simeq \hat{e}\circ h$,then $g
\simeq h$.

By the result of Sullivan, it follows immediately that the map
\[
\hat{e}_*:map_*(X^n,Y)_{f|_{X^n}} \to
map_*(X^n,\hat{Y})_{\hat{f}|_{X^n}}
\]
has also the property above. Thus the induced homomorphism
\[
(\hat{e}_*)_{\#}:\pi_j(map_*(X^n,Y);f|_{X^n}) \to
\pi_j(map_*(X^n,\hat{Y});\hat{f}|_{X^n})
\]
is injective. This completes the proof.

Next we show that if  $(\tau^*)_{\#}(\alpha)=0$,then  $\alpha \in
\pi_j(map_*(X,Y);f)$ is a phantom element. From the assumption, we
have
\[
(\hat{e}_*)_{\#}(\tau^*)_{\#}( \alpha)=0 \in
\pi_j(map_*(X_{tau},\hat{Y});\hat{f}_{\tau})
\]
By Proposition2.1 of \cite{os}, $map_*(X_{(0)},\hat{Y})$ is weakly
contractible and hence the natural map
$\pi_j(map_*(X,\hat{Y});\hat{f}) \to
\pi_j(map_*(X_{\tau},\hat{Y});\hat{f}_{\tau})$ is an isomorphism
for $j > 0$. Since $(\hat{e}_*)_{\#}\circ
(\tau^*)_{\#}=(\tau^*)_{\#} \circ (\hat{e}_*)_{\#} $,we have
\[
\hat{\alpha}=0 \in \pi_j(map_*(X,\hat{Y});\hat{f})
\]
By the first part, $\alpha$ is a phantom element.
\end{proof}

\begin{prop}\label{T:fprop}
Let $X,Y$ be CW complexes such that $[\Sigma^jX_{\tau},Y]=0$ and
$[\Sigma^{j+1}X_{\tau},Y]=0$. If $f:X \to Y$ is a phantom map,then
we have
\[
 Ph^f_j(X,Y)= \pi_j(map_*(X,Y);f)=\pi_j(map_*(X_{(0)},Y);f_{(0)})
\]
\end{prop}
\begin{proof}
To show the first equation, it suffices to prove by
Theorem\ref{T:fundamental} that
\[
\pi_j(map_*(X_{\tau},Y);f_{\tau})=0
\]
 From Theorem 5.1 of
\cite{cm95} we have $f_{\tau}\simeq0$ and this implies
\[
\pi_j(map_*(X_{\tau},Y);f_{\tau})\cong
\pi_j(map_*(X_{\tau},Y);*)\cong [\Sigma^jX_{\tau},Y]=0
\]
which is what we want to prove.

Similarly we can prove $\pi_{j+1}(map_*(X_{\tau},Y);f_{\tau})=0$ .
By using the fibration
\[
map_*(X_{(0)},Y)\to map_*(X,Y) \to map_*(X_{\tau},Y)
\]
and the fact that
\[
\pi_j(map_*(X_{\tau},Y);f_{\tau})=\pi_{j+1}(map_*(X_{\tau},Y);f_{\tau})=0
\]
we obtain immediately the second equation.
\end{proof}
\begin{cor}\label{T:liegroup}
Let $X,Y,f$ be as in Proposition\ref{T:fprop}. If we assume
further that $Y$ is a 1-connected rational $H$-space. Then
\[
Ph^f_j(X,Y)=\pi_j(map_*(X,Y);f)=[\Sigma^jX_{(0)},Y]
\]
\end{cor}
Before the proof of the Corollary  let us first recall some useful
results.

\begin{lem}\cite{az}
Let $X,Y_1,Y_2$ be 1-connected CW complexes with finite type and
$t:Y_1 \to Y_2$ be a rational equivalence. Then the induced map
$t_*: map_*(X_{(0)},Y_1) \to map_*(X_{(0)},Y_2)$ is a weak
equivalence.
\end{lem}
\begin{proof}[Proof of Corollary\ref{T:liegroup}]
Let $t:Y \to \bar{Y}$ be an integral approximation of $Y$(see
\cite{az} for the concept of integral approximation of a
space).Then by the above Lemma , the map $t_*: map_*(X_{(0)},Y)
\to map_*(X_{(0)},\bar{Y})$ is a weak equivalence. Since $Y$ is a
rational H-space,we can choose $\bar{Y}$ to be an H-space(see
remark by McGibbon in\cite{cm95}). It is well known that the
different components of a mapping space $ map_*(X_{(0)},\bar{Y})$
have the same homotopy type if $\bar{Y}$ is an H-space. This
completes the proof.
\end{proof}
\begin{thm}\label{T:section2maincoro1}
Let $X=K(H,m+1)$ , $Y=Baut(P)$ and $f:X \to Y $ is any map where
$P$ is a simply connected finite CW complex and  $m \geq 2$.Then
\[
 Ph^f_j(X,Y)= \pi_j(map_*(X,Y);f)=
 =[\Sigma^jX_{(0)},Y]
\]
\end{thm}
\begin{proof}
 We know from Corollary C' of \cite{az} that $map_*(X_{\tau},Y)$ is
weakly contractible.It follows that any map $f:X \to Y $ is a
phantom map. Thus we can apply the Proposition\ref{T:fprop} to get
\[
Ph_j^f(X,Y)=\pi_j(map_*(X,Y);f)=\pi_j(map_*(X_{(0)},Y);f_{(0)})
\]
To prove the last equation, note that, if $X$ is 1-connected, then
we have
\[
map_*(X,Y)\simeq map_*(X,\tilde{Y})
\]
where $\tilde{Y}$ is the universal covering of $Y$. Now $Y$ and
thus $\tilde{Y}$ is a 1-connected rational H-space.It follows from
Corollary\ref{T:liegroup}  and the above equivalence that
\[
\pi_j(map_*(X_{(0)},Y);f_{(0)})=\pi_j(map_*(X_{(0)},Y);*)
\]
\end{proof}
\begin{thm}\label{T:liephan}
 Let $X=BG$, $Y=Baut(P)$ and $f:X \to Y$ is a phantom map where $G$ is a connected compact Lie
group and $P$ is 1-connected finite CW complex. Then for $j \geq
1$ we have
\[
Ph_j^f(X,Y)=\pi_j(map_*(X,Y);f)= [\Sigma^jX_{(0)},Y]
\]
\end{thm}
\begin{proof}
By Proposition\ref{T:fprop},to prove the first equation it
suffices to prove that
$[\Sigma^jX_{\tau},Y]=[\Sigma{j+1}X_{\tau},Y]=0$. Now for $i \geq
1$, we have
\[
[\Sigma^iBG_{\tau},Y]=[\Sigma^{i-1}BG_{\tau},\Omega Y]=
\]
\[
=[\Sigma^{i-1}BG_{\tau},aut(P)]=0
\]
where the last equation follows from Theorem C(c) of
Zabrodsky\cite{az}.

The proof of the second equation follows from the same argument as
in the proof of Theorem\ref{T:section2maincoro1}.
\end{proof}

\section{Forgetable map and its description}
Let us consider the principal G- bundle $q:P\to B$ with structure
group G where G acts on P freely. For each such bundle one can
consider the space $aut^{G}(P)$ of unbased G-equivariant
self-homotopy equivalences of P and the group
\[
Aut^{G}(P)=\pi_{0}(aut^{G}(P))
\]
which is called the group of G-equivariant self-homotopy
equivalences of the principal bundle $q:P\to B$. On the other
hand, we can also consider the space aut(P) of unbased
self-homotopy equivalences of space P and the group
\[
Aut(P)=\pi_{0}(aut(P))
\]
which is called the group of unbased self-homotopy equivalences of
P. There have been extensive study on these two subjects , see
\cite{dwk} and the extensive references there.About ten years ago,
Tsukiyama \cite{dwk} posed the following
\begin{ques}
When is the natural map $F:Aut^{G}(P) \to Aut(P)$, which forgets
the G-action, a monomorphism?
\end{ques}

No progress has been made unless Tsukiyama's two recent papers
\cite{kt1} \cite{kt2}. In this paper, we will try to attack this
question . Our approach is based on the identification of the
space of G-equivariant self-homotopy equivalences as the loop
space on a mapping space and the recent results on the Sullivan
conjecture. What is most interesting about our results is the
relation between the injectivity of F and the existence of the
phantom map between appropriate spaces.

In \cite{kt1} \cite{kt2}, Tsukiyama used an indirect approach to
attack the question and got partial results on it.In this section,
based on a simple but crucial observation, we will identify the
homomorphism F as the homomorphism induced on $\pi_0$ by a map
whose homotopy fiber can be determined explicitly and thus can
determine the kernel of F under reasonable condition. Now let G be
a topological group , $q:P\to B$ be a principal G-bundle and $k:B
\to BG$ be the its classifying map. For the map k ,we can take
$\bar{k}:\bar{B} \to BG$ as a fibration via the standard
factorization of a map into the composite of a homotopy
equivalence and a Hurewicz fibration. Given fibration
$\bar{k}:\bar{B} \to BG$, we can form the group
$Aut_{BG}(\bar{B})$ the group of homotopy classes of self homotopy
equivalences of $\bar{B}$ over $BG$.The following is a  well known
result  \cite{bhmp},\cite{gottlieb2},\cite{gottlieb3}.
\begin{prop}
There is a natural isomorphism
\[
Aut^{G}(P) \cong Aut_{BG}(\bar{B}) \cong \pi_1(Map(BG,Baut(P)),c)
\]
where  map $c:BG \to Baut(P)$ is determined by the principal
bundle.
\end{prop}
If the above isomorphism is natural in object G ,then  the map $F$
will be naturally isomorphic to the map
\[
ev_*:\pi_1(Map(BG,Baut(P)),c) \to \pi_1Baut(P)
\]
whose kernel can be computed explicitly by the evaluation
fibration
\[
Map_*(BG,Baut(P))_c \to Map(BG,Baut(P))_c \to Baut(P)
\]
A careful check confirms the above speculation and leads to the
following which is the Theorem\ref{T:main5} in the Introduction.
\begin{thm} \label{T:main}
Let $q: P \to B $ be a principal G-bundle, $c:BG \to Baut(P)$ be
the classifying map for fibration $\bar{k}:\bar{B} \to BG$. Then
there is a commutative diagram
\[
\begin{CD}
\pi_0aut_{BG}(\bar{B)} @>>\cong> \pi_0aut^G(P)\\
 @VVV    @VVV \\
\pi_0aut(\bar{P}) @>>\cong>\pi_0aut(P)
\end{CD}
\]
where $\bar{P}$ is the fiber of the fibration $\bar{k}$ which is
homotopy equivalent to $P$ and the two horizontal maps are
isomorphisms. It follows from the diagram above immediately that
the following sequence is exact
\[
\pi_1aut(P) \to \pi_1(Map_*(BG,Baut(P)),c) \to Aut^G(P)) \to
Aut(P)
\]
\end{thm}
This theorem follows directly from the following lemmas.
 Let $q: P \to B$ be a  principal G-bundle
and $k: X \to B$ be a map, then there is an associated  principal
G-bundle over X defined by
\[
\begin{CD}
k^*(P) @>>> P\\
 @VVk^*(q)V    @VVqV \\
X @>>k> B
\end{CD}
\]

\begin{lem}\label{T:lemma1}
Let $\pi: EG \to BG$ be the universal principal G-bundle,then the
rule that takes $k$ to $k^*(\pi)$ defines a natural bijection from
$[B,BG]$, the set of free homotopy classes of maps from $B$ to
$BG$, to the set of isomorphism classes of G-bundles over B.
\end{lem}
\begin{proof}
well known.
\end{proof}

\begin{lem} \label{T:lemma2}
If $q: P \to B$ is a principal G-bundle and $g: X \to B$ is a
homotopy Equivalence, then the induced bundle map from $
g^*(q):{g^{*}(P)} \to X$ to q is a homotopy equivalence between
two principal bundles.
\end{lem}

\begin{proof}
This is (1.9)of \cite{pb}.
\end{proof}

\begin{lem} \label{T:lemma3}
Let $q: P \to B$ and  $k:B \to BG$ be as above.Then there is a
commutative diagram  .
 \[
\begin{CD}
aut^G(k^*(EG)) @>>> aut^G(P)\\
 @VVV    @VVV \\
 aut(k^*(EG)) @>>> aut(P)
\end{CD}
\]
Where the two vertical maps are Forgetable maps and the horizontal
maps are homotopy equivalences and are defined in the proof.
\end{lem}

\begin{proof}
 By Lemma \ref{T:lemma1}, There is a principal bundle isomorphism over B
$h:P \to k^*(EG)$ . Define the horizontal maps by the rule
\[
x \longmapsto h^{-1}\circ x \circ h
\]
It is obvious that the diagram is commutative.
\end{proof}

\begin{lem} \label{T:lemma4}
Let $q: P \to B$ , $k:B \to BG$ and $\bar{k}:\bar{B} \to BG$ be as
above.Then there is a commutative diagram up to homotopy .
 \[
\begin{CD}
aut^G(\bar{k}^*(EG)) @>>> aut^G(k^*(EG)) \\
 @VVV    @VVV \\
 aut(\bar{k}^*(EG)) @>>> aut(k^*(EG))
\end{CD}
\]
Where the two vertical maps are Forgetable maps and the horizontal
maps are homotopy equivalences and are defined in the proof.
\end{lem}

\begin{proof}
By Lemma \ref{T:lemma2} ,There is a homotopy equivalence h of two
principal bundles $k^*(\pi):k^*(EG) \to B $ and
$\bar{k}^*(\pi):\bar{k}^*(EG) \to \bar{B}$. As in the proof of the
lemma above , define the horizontal maps similarly. Then the
diagram is easily seen to be commutative up to homotopy.

\end{proof}

\begin{lem} \label{T:lemma5}
Let $q: P \to B$ , $k:B \to BG$ and $\bar{k}:\bar{B} \to BG$ be as
above.Then there is a commutative diagram up to homotopy .
 \[
\begin{CD}\label{T:diagram}
\pi_0aut_{BG}(\bar{B}) @>>> \pi_0aut^G(\bar{k}^*(EG))
\\
 @VVV    @VVV \\
 \pi_0aut(\bar{P}) @>>> \pi_0aut(\bar{k}^*(EG))
\end{CD}
\]
Where the $\bar{P}$ is the fiber of the fibration $\bar{k}$ which
is homotopy equivalent to P,the right vertical map is forgetting
maps, the left
 vertical map is the map by taking a $f$ to the map which induced
on the fiber of the fibration $\bar{k}$ at a based point of $BG$
and the horizontal maps are defined in the proof.
\end{lem}

\begin{proof}
Consider the following diagram
 \[
\begin{diagram}
\node[2]{\bar{P}} \arrow{e,t}{f} \arrow{ssw,l}{\tilde{h}}
\arrow{se} \node{\bar{k}^*(EG)} \arrow[2]{e} \arrow{s}
\arrow{ssw,l}{\bar{h}} \node[2]{EG} \arrow{s} \arrow{ssw,r}{=}\\
\node[3]{\bar{B}} \arrow[2]{e} \arrow{ssw,l}{h} \node[2]{BG}
\arrow{ssw,r}{=}\\ \node{\bar{P}} \arrow{e,t}{f} \arrow{se}
\node{\bar{k}^*(EG)} \arrow{s} \arrow[2]{e} \node[2]{EG}
\arrow{s}\\ \node[2]{\bar{B}} \arrow[2]{e} \node[2]{BG}
\end{diagram}
\]
By definition
\[
\bar{k}^*(EG)=\{(b,e)\in \bar{B}\times EG |\bar{k}(b)=\pi(e)\}
\]

\[
\bar{P}=\{(b,*)\in \bar{B}\times EG |\bar{k}(b)=*\}
\]
Now there is an obvious map
\[
f:\bar{P} \to \bar{k}^*(EG), (b,*)\longmapsto (b,*)
\]
which is a homotopy equivalence by the general property of
pullback.

Given a self homotopy equivalence $h \in aut_{BG}(\bar{B})$ there
exists a map
\[
\bar{h}:\bar{k}^*(EG)\to \bar{k}^*(EG)
\]
defined by $\bar{h}(b,e)=(h(b),e)$ which makes the pair
$(\bar{h},h)$ a principal bundle map. The given map $h$ induces
also an obvious map $\tilde{h}\in aut(\bar{P})$ defined by
\[
\tilde{h}(b,*)=(h(b),*)
\]

It is easy to check that $\bar{h}\circ f=f \circ \tilde{h}$

If we define the horizontal maps in the diagram by
\[
\pi_0aut_{BG}(\bar{B}) \to \pi_0aut^G(\bar{k}^*(EG)), h
\longmapsto (\bar{h},h)
\]
\[
\pi_0aut(\bar{P}) \to \pi_0aut(\bar{k}^*(EG)), g \longmapsto f
\circ g \circ f^{-1}
\]
then it is easy to check that the diagram is commutative which is
what want to prove .
\end{proof}

\section{Applications  to the problem of Forgetable maps}
In the last section we have embedded the Forgetable map into the
exact sequence
\[
\pi_1aut(P) \to \pi_1(Map_*(BG,Baut(P)),c) \to Aut^G(P)) \to
Aut(P)
\]

In this section we will apply the phantom map theory to extract
information about the Forgetable map.If

\[
\pi_1(Map_*(BG,Baut(P)),c)=0
\]
then we know that  $KerF=0$ from the above exact sequence. If
\[
\pi_1(Map_*(BG,Baut(P)),c)\neq0
\]
then we can't say anything about the $KerF$. On the another hand
$\pi_1aut(P)$ is a countable group if $P$ is a finite CW complex.
It follows that $KerF$ is uncountable if
\[
\pi_1(Map_*(BG,Baut(P)),c)
\]
is uncountable. This is the point where we find the relation
between phantom map theory and Forgetable map. From the discussion
above , we have the following which is Theorem\ref{T:4} in the
Introduction.
\begin{thm}\label{T:3.1}
Let $\pi:P \to B$ be  as above , $k:B \to BG$ the classifying map,
$\bar{k}:\bar{B} \to BG$ the associated fibration with fiber
$\bar{P}$ and $c:BG \to Baut(P)$ the classifying map of $\bar{k}$.
The following are true.
\begin{itemize}
 \item{If $KerF$ is uncountable if $Ph_1^c(BG,Baut(P))\neq0$}
 \item{If $c$ is phantom map and $G$ be a connected compact Lie group , $KerF$ is either zero or uncountable}
\item{If $G=K(H,m)$ where $m \geq 2$ and $H$ is a finitely generated abelian group,$KerF$ is either zero or uncountable}
\end{itemize}
\end{thm}
\begin{proof}
It suffices to prove the last two statements.

By Theorem\ref{T:section2maincoro1} when $G$ be a 1-connected
K(H,m) or Theorem\ref{T:liephan} when $c$ is phantom map and $G$
be a connected compact Lie group , we have
\[
Ph^c_1(BG,Baut(P))=\pi_1(map_*(BG,Baut(P));c)
\]
Thus $\pi_1(map_*(BG,Baut(P));c)$ is either zero or uncountable
\end{proof}
\begin{cor}\label{T:cor**}
If the map $c$   in the Theorem \ref{T:3.1} is a phantom map  ,
then $KerF$ is uncountable iff $[BG_{(0)},aut(P)]$ is nontrivial.
\end{cor}
\begin{proof}
By the Theorem\ref{T:liephan} for $j=1$, we have
\[
Ph_j^c(BG,Baut(P))=\pi_j(map_*(BG,Baut(P));c)
\]
Thus the kernel of the Forgetable map is either zero or
uncountable. It is zero iff
\[
\pi_j(map_*(BG,Baut(P));c)=0
\]
By the same Theorem again for $j=1$, we have
\[
\pi_1(map_*(BG,Baut(P));c)=[\Sigma^1BG_{(0)},Baut(P)]=[BG_{(0)},aut(P)]
\]
It follows that kernel of Forgetable map is zero iff
\[
[BG_{(0)},aut(P)]=0
\]
\end{proof}
The following corollary is   Theorem\ref{T:5} in the Introduction.
\begin{cor}\label{T:cor*}
Let $P$ be a 1-connected CW complex  . Then there is  a connected
compact Lie group $G$ and a principal $G$-bundle  such that the
total space has the homotopy type of $P$ , the classifying map c
is a phantom map and the associated Forgetable map $F$ has
uncountable kernel iff
$\bigoplus_{i>0}\pi_{2i+1}(map(P_{(0)},P_{(0)});id)$ is
nontrivial.
\end{cor}
\begin{proof}
"$\Leftarrow$":First note that given map $c_0:BG \to Baut(P)$
there is a principal bundle such that the total space has the
homotopy type of $P$ and the natural associated map $c$ is
homotopy to the given $c_0$.

To prove the Corollary it is sufficient to take $c_0=*$ and choose
a Lie group G such that
$[BG_{(0)},aut(P)]=[BG_{(0)},map(P,P)_{id}]\neq 0$.

According to Theorem\ref{T:thm1.4},we have  for some even integer
$t>0$,
$[BG_{(0)},map(P,P)_{id}]\neq 0$ if $H^t(BG,Q)\neq 0$ and
$\pi_{t+1}(map(P_{(0)},P_{(0)});id)\neq0 $

Let $t_0$ be the smallest positive even integer such that
\[
\pi_{t+1}(map(P_{(0)},P_{(0)});id)\neq0
\]
 There exists of course a
compact Lie group $G$ such that $H^{t_0}(BG,Q)\neq 0$.

It follows from the discussion above that there exists principal
$G$-bundle such that the total space has the homotopy type of $P$
and the associated Forgetable map $F$ has uncountable kernel.

"$\Rightarrow$":If
$\bigoplus_{i>0}\pi_{2i+1}(map(P_{(0)},P_{(0)});id) =0$, it is
easy to see that $[BG_{(0)},aut(P)]=0$ for any connected compact
Lie group. This completes the proof by the Corollary\ref{T:cor**}

\begin{cor}
Let $P$ be a 1-connected finite CW complex.Then for all principle
$G$-bundles $q:P \to B$ such that the structure group is a
connected compact Lie group and the associated map $c$ is
phantom,the associated Forgetable map is injective iff
$\bigoplus_{i>0}\pi_{2i+1}(map_*(P_{(0)},P_{(0)});id)=0$
\end{cor}
\end{proof}

Similarly to th Corollary\ref{T:cor*} we  have the following when
$G=K(H,m)$ .
\begin{cor}\label{T:cor***}
Let $P$ be a 1-connected CW complex  . Then there is some
principal $K(H,2m+1)$(respectively,$K(H,2m)$)-bundle,$m \geq 1$,
such that the total space has the homotopy type of $P$  and the
associated Forgetable map $F$ has uncountable kernel iff
$\bigoplus_{i>1}\pi_{2i+1}(map(P_{(0)},P_{(0)});id)$ (respectively
, $\bigoplus_{i>1}\pi_{2i}(map(P_{(0)},P_{(0)});id)$) is
nontrivial.
\end{cor}
\begin{cor}
Let $P$ be a 1-connected CW complex . Then for all $m \geq 1$ ,
finitely generated abelian group $H$ and every principal
$K(H,2m+1)$ -bundle with total space homotopy equivalent to $P$
,the associated Forgetable map is injective iff
$\bigoplus_{i>1}\pi_{2i+1}(map(P_{(0)},P_{(0)});id)=0$ .
\end{cor}
The following Corollary is the Theorem\ref{T:1} in the
Introduction

\begin{cor}
Let $P$ be a 1-connected CW complex . Then for all $m \geq 1$ ,
finitely generated abelian group $H$ and every principal
$K(H,2m)$-bundle with total space homotopy equivalent to $P$ ,the
associated Forgetable map is injective iff
$\bigoplus_{i>1}\pi_{2i}(map(P_{(0)},P_{(0)});id)=0$ .
\end{cor}

Let us conclude this paper with another question motivated by the
results obtained in this paper.
\begin{ques}
Is it possible that for every 1-connected CW complex $P$ there
exist a compact Lie group $G$ and a principal $G$-bundle such that
the total space has the homotopy type of $P$ and the associated
Forgetable map $F$ has uncountable kernel?
\end{ques}

----------------------------------------------------------------

----------------------------------------------------------------

----------------------------------------------------------------

\end{document}